\documentclass[10pt]{article}
\usepackage{a4}
\usepackage[margin=0.82in]{geometry}
\usepackage{amsthm, amsmath, amsfonts, cite}
\usepackage{graphicx}
\usepackage[margin=0.8cm]{caption}
\usepackage{titlesec}
\usepackage{epsfig}
\usepackage{caption}  
\usepackage{capt-of}
\usepackage{tikz}
\usepackage{multicol}
\usepackage{thmtools}
\usepackage{thm-restate}
\usepackage{cleveref}
\usepackage{enumerate}
\usepackage{braket}
\newtheorem{theorem}{Theorem}
\newtheorem{proposition}[theorem]{Proposition}
\newtheorem{lemma}[theorem]{Lemma}

\newtheorem{claim}{Claim}[theorem]

\allowdisplaybreaks

\title{Odd distances in colourings of the plane}
\author{James Davies\thanks{Department of Pure Mathematics and Mathematical Statistics, Wilberforce Road, Cambridge CB3 0WB, UK. E-mail: \texttt{jgd37@cam.ac.uk}.}}
\date{}

\begin{document}

\maketitle

\begin{abstract}
	We prove that every finite colouring of the plane contains a monochromatic pair of points at an odd distance from each other.
\end{abstract}

\section{Introduction}

The \emph{odd distance graph} is the graph with vertex set $\mathbb{R}^2$, and whose edges are the pairs of points in $\mathbb{R}^2$ whose Euclidean distance is an odd integer.
The \emph{chromatic number} of a graph is equal to the minimum number of colours required to assign each vertex a colour so that no two adjacent vertices receive the same colour.
In 1994, Rosenfeld posed to Erd\H{o}s the problem of determining the chromatic number of the odd distance graph, and Erd\H{o}s~\cite{erdos1994twenty} presented the problem the next day in an open problems talk. 
Since then, the problem has been reiterated numerous times (see for instance,~\cite{ardal2009odd,bachoc2014spectral,beineke2015topics,chybowska2022coloring,damasdi2021odd,golovanov2023odd,heuleeasier,jensen2011graph,kechris2016descriptive,parts20226,quas2009distances,rosenfeld1996odd,rosenfeld2008some,rosenfeld2014forbidden,soifer2010between,soifer2016hadwiger,steinhardt2009coloring,kalai15some,heuleodd,graham2010open}).

At the time the best known lower bound was $4$, the same as for the Hadwiger-Nelson problem. In 2009, Ardal, Ma{\v{n}}uch, Rosenfeld, Shelah, and Stacho~\cite{ardal2009odd} constructed a 5-chromatic subgraph of the odd distance graph. Making a breakthrough in the Hadwiger-Nelson problem, de Grey~\cite{de2018chromatic} further constructed a 5-chromatic unit distance graph.
Motivated by work towards constructing a 6-chromatic unit distance graph, Heule~\cite{heuleeasier,heuleodd} recently offered \$500 for determining whether or not the odd distance graph is 5-colourable. The prize was promptly claimed by Parts~\cite{parts20226} with a construction of a 6-chromatic subgraph of the odd distance graph.

As conjectured by Soifer~\cite{soifer2016hadwiger}, we prove that the odd distance graph has no finite colouring.

\begin{theorem}\label{odd}
	For every finite colouring of the plane, there is a monochromatic pair of points at an odd distance from each other.
\end{theorem}

Since there are countable colourings of the odd distance graph, the chromatic number of the odd distance graph is therefore countably infinite. This solves Rosenfeld's odd distance problem. An example of such a countable colouring is given by assigning the points $[\frac{1}{2}x, \frac{1}{2}(x+1)) \times [\frac{1}{2}y, \frac{1}{2}(y+1)) \subset \mathbb{R}^2$ the colour $(x,y)$ for each $(x,y)\in \mathbb{Z}^2$.

In fact, we prove the following more general result.

\begin{theorem}\label{p mod q}
	For every finite colouring of the plane, for every pair of integers $p,q$, there is a monochromatic pair of points at a distance congruent to $p \pmod{q}$ from each other.
\end{theorem}

Very recently, with Rose McCarty and Micha\l{} Pilipczuk~\cite{davies2023prime} we have extended the methods in the present paper to prove two generalizations of Theorem~\ref{odd} for prime and polynomial distances.

Our construction may also be of independent interest from a purely graph theoretic point of view.
Indeed, when $q$ does not divide $p$, the graphs we construct to prove Theorem~\ref{p mod q} happen to be triangle-free and appear to not contain any of the classical constructions of triangle-free graphs with large chromatic number (see~\cite{scott2020survey} for a discussion of such constructions).

If we restrict our colourings of the plane to be measurable, then much more than Theorem~\ref{p mod q} is known.
F{\"u}rstenberg, Katznelson, and Weiss~\cite{furstenberg1990ergodic} gave an ergodic-theoretic proof that for every finite measurable colouring of the plane, there exists a $d_0>0$ such that for every real $d\ge d_0$, there is a monochromatic pair of points at a distance of $d$ from each other.
Later, Bourgain~\cite{bourgain1986szemeredi} gave a harmonic-analytic proof, Falconer and Marstrand~\cite{falconer1986plane} gave a more direct geometric proof, and Oliveira and Vallentin~\cite{de2010fourier} gave a proof using Fourier analysis.

Steinhardt~\cite{steinhardt2009coloring} gave a spectral proof of Theorem~\ref{odd} under the additional assumption of the colouring being measurable.
Inspired by this, our proof of Theorem~\ref{p mod q} also uses spectral techniques. In Section~2 we prove an analogue of the Lov{\'a}sz theta bound~\cite{lovasz1979shannon} for Cayley graphs of $\mathbb{Z}^d$ (see Theorem~\ref{Z^d theta}).
Then in Section~3, we prove Theorem \ref{p mod q} by constructing Cayley graphs of $\mathbb{Z}^2$ (see Theorem~\ref{graphs}) which can both be realized as distance graphs in the plane, and for which we can effectively apply Theorem~\ref{Z^d theta}.

\section{Ratio bound}

We say that a set $C\subseteq \mathbb{Z}^d$ is \emph{centrally symmetric} if $C=-C$. Similarly, a function $w: C \to \mathbb{R}$ is \emph{centrally symmetric} if $w(x)=w(-x)$ for all $x\in C$. For a centrally symmetric set $C\subseteq \mathbb{Z}^d \backslash \{0\}$, we let $G(\mathbb{Z}^d, C)$ be the Cayley graph of $\mathbb{Z}^d$ with generating set $C$.
This is the graph with vertex set $\mathbb{Z}^d$, where two vertices $u,v\in \mathbb{Z}^d$ are adjacent if there exists a $c\in C$ with $u+c=v$.
Given a real symmetric matrix $B$, we let $\lambda_{\max}(B)$ and $\lambda_{\min}(B)$ be equal to its maximum and minimum eigenvalue respectively.
For $x,y\in \mathbb{Z}^d$, we write $x \cdot y$ for the dot product.

An \emph{independent set} of a graph $G$ is a set of pairwise non-adjacent vertices. The \emph{independence number} of a graph $G$ is equal to the maximum size of an independent set in $G$.
For a graph $G$, the independence number and chromatic number of $G$ are denoted by $\alpha(G)$ and $\chi(G)$ respectively.
For a finite graph $G$, we let $\overline{\alpha}(G) = \alpha(G)/|V(G)|$.
Notice that $\chi(G) \ge 1/ \overline{\alpha}(G)$. For Cayley graphs on $\mathbb{Z}^d$, we define $\overline{\alpha}(G)$ in terms of the maximum upper density of an independent set.
Given a set $I\subseteq \mathbb{Z}^d$, its \emph{upper density} is
\[
\delta(I) = \limsup_{R \to \infty} \frac{|I \cap [-R,R]^d  |}{(2R+1)^d}.
\]
For a Cayley graph $G(\mathbb{Z}^d , C)$, we let
\[
\overline{\alpha}(G(\mathbb{Z}^d, C)) = \sup \{ {\delta(I)} : I \text{ is an independent set of } G(\mathbb{Z}^d, C) \}.
\]
If $V(G(\mathbb{Z}^d, C))= \mathbb{Z}^d$ can be partitioned into $k$ independent sets $I_1, \ldots, I_k$, then $k\overline{\alpha}(G) \geq \delta(I_1)+\dots+\delta(I_k) \geq \delta(V(G))=1$.
In particular, notice that $\chi(G(\mathbb{Z}^d, C)) \ge 1/ \overline{\alpha}(G(\mathbb{Z}^d, C))$.

The purpose of this section is to prove the following analogue of the Lov{\'a}sz theta bound~\cite{lovasz1979shannon}.

\begin{theorem}\label{Z^d theta}
	Let $C\subseteq \mathbb{Z}^d \backslash \{0\}$ and $w:C \to \mathbb{R}$ be centrally symmetric with $\sum_{x\in C} w(x)$ absolutely convergent and positive. Then
	\[
	\overline{\alpha}(G(\mathbb{Z}^d, C)) \le \frac{-\inf_{u\in (\mathbb{R} / \mathbb{Z})^d}  \widehat{w}(u) }{  \sup_{u\in (\mathbb{R} / \mathbb{Z})^d}  \widehat{w}(u)   -\inf_{u\in (\mathbb{R} / \mathbb{Z})^d}  \widehat{w}(u)},
	\]
	where
	\[
	\widehat{w}(u) = \sum_{x\in C} w(x) e^{2\pi i u \cdot x}.
	\]
\end{theorem}

We remark that one could derive Theorem~\ref{Z^d theta} from a measurable analogue for Cayley graphs $G(\mathbb{R}^d,C)$ due to Bachoc, DeCorte, Oliveira, and Vallentin~\cite{bachoc2014spectral}.
Indeed, Dutour Sikiri{\'c}, Madore, Moustrou, and Vallentin~\cite{dutour2021coloring} use this to derive a slightly less general version of Theorem~\ref{Z^d theta} which gives a bound for the chromatic number of the Voronoi tessellation of a lattice.
For ease of accessibility, we will instead prove Theorem~\ref{Z^d theta} directly from Lov{\'a}sz's~\cite{lovasz1979shannon} theta bound for finite graphs, and a simple lemma computing the eigenvalues of a weighted adjacency matrix of a Cayley graph $G(\mathbb{Z}^d_n , C)$, where $\mathbb{Z}^d_n$ is the abelian group $( \mathbb{Z}/n\mathbb{Z} )^d$.

\begin{theorem}[Lov{\'a}sz~\cite{lovasz1979shannon}]\label{finite theta}
	If $G$ is a finite graph, then
	\[
	\overline{\alpha}(G) \le \min_B \frac{-\lambda_{\min}(B)}{\lambda_{\max}(B) - \lambda_{\min}(B)},
	\]
	where $B$ ranges over all symmetric matrices in $\mathbb{R}^{V(G)\times V(G)}$ such that $B_{u,v} = 0$ whenever $uv\not\in E(G)$.
\end{theorem}

\begin{lemma}[Lov{\'a}sz~\cite{lovasz1975spectra}]\label{eigenvalues}
	Let $w:\mathbb{Z}_n^d \to \mathbb{R}$ be centrally symmetric, and let $B$ be the symmetric matrix in $\mathbb{R}^{\mathbb{Z}_n^d \times \mathbb{Z}_n^d}$ with $B_{u,v}= w(u-v)$ for each pair $u,v\in \mathbb{Z}_n ^d$. Then the eigenvalues of $B$ are $\{\lambda_z : z\in \mathbb{Z}_n^d \}$, where for each $z\in \mathbb{Z}_n^d$, 
	\[
	\lambda_z = \sum_{x\in \mathbb{Z}_n^d} w(x)e^{\frac{2\pi i}{n} z\cdot x}.
	\]
\end{lemma}

We now derive Theorem~\ref{Z^d theta} from Theorem~\ref{finite theta} and Lemma~\ref{eigenvalues}.

\begin{proof}[Proof of Thereom~\ref{Z^d theta}.]
	Since $\sum_{x\in C} w(x)$ is absolutely convergent and positive, we have $\sum_{x\in C} |w(x)| = s$ for some finite real $s>0$.
	For each real $0 < \epsilon < s$, let $C_\epsilon$ be a finite subset of $C$ such that $\sum_{x\in C_\epsilon} |w(x)| > s - \epsilon$.
	For $n$ sufficiently large,
	let $w_\epsilon:\mathbb{Z}^d_n \to \mathbb{R}$ be such that $w_\epsilon(x)=0$ for $x\not\in C_\epsilon$, and $w_\epsilon(x)=w(x)$ for $x\in C_\epsilon$, and let $B_{\epsilon,n}$ be the symmetric matrix in $\mathbb{R}^{\mathbb{Z}_n^d \times \mathbb{Z}_n^d}$ with $B_{u,v}= w_\epsilon(u-v)$ for each pair $u,v\in \mathbb{Z}_n ^d$.
	For $u\in (\mathbb{R} / \mathbb{Z} )^d$, let
	\[
	\widehat{w}_\epsilon(u) = \sum_{x\in C_\epsilon} w_\epsilon(x) e^{2\pi i u \cdot x}.
	\]
	By Lemma~\ref{eigenvalues}, the eigenvalues of $B_{\epsilon,n}$ are $\{\lambda_z : z\in \mathbb{Z}_n^d \}$, where for each $z\in \mathbb{Z}_n^d$, 
	\[
	\lambda_z = \sum_{x\in \mathbb{Z}_n^d} w_\epsilon(x)e^{\frac{2\pi i}{n} z\cdot x}= 
	\sum_{x\in C_\epsilon} w_\epsilon(x)e^{\frac{2\pi i}{n} z\cdot x} = \widehat{w}_\epsilon(z/n).
	\]
	So, it follows that $\lim_{n \to \infty} \lambda_{\max}(B_{\epsilon,n}) = \sup_{u\in (\mathbb{R} / \mathbb{Z})^d}  \widehat{w}_\epsilon(u)$, and $\lim_{n \to \infty} \lambda_{\min}(B_{\epsilon,n} ) = \inf_{u\in (\mathbb{R} / \mathbb{Z})^d}  \widehat{w}_\epsilon(u)$.
	In particular, by Theorem~\ref{finite theta}, for each real $0 < \epsilon < s$, we have that \[
	\limsup_{n \to \infty} \overline{\alpha}(G(\mathbb{Z}_n^d ,C_\epsilon ))
	\le
	\limsup_{n \to \infty} \frac{-\lambda_{\min}(B_{\epsilon,n})   }{  \lambda_{\max}(B_{\epsilon,n})  -\lambda_{\min}(B_{\epsilon, n})}
	=
	\frac{-\inf_{u\in (\mathbb{R} / \mathbb{Z})^d}  \widehat{w}_\epsilon(u) }{  \sup_{u\in (\mathbb{R} / \mathbb{Z})^d}  \widehat{w}_\epsilon(u) -\inf_{u\in (\mathbb{R} / \mathbb{Z})^d}  \widehat{w}_\epsilon(u)}.\]
	
	Consider some real $0 < \epsilon < s$. Then, since $|C_\epsilon|$ is finite, there exists some positive integer $T$ such that $C_\epsilon \subset [-T,T]^d$.
	For every real $\epsilon' >0$ and positive integer $N$, there exists an independent set $I$ of $G(\mathbb{Z}^d ,C_\epsilon )$ and a positive integer $n\ge N$ with $\frac{|I \cap [-n,n]^d|}{(2n+1)^d} \ge \overline{\alpha}(G(\mathbb{Z}^d ,C_\epsilon )) - \epsilon'$.
	Then, $I_n = I \cap [-n,n]^d$ is an independent set of $G(\mathbb{Z}^d_{2n+1 + T},C_\epsilon )$ since the subgraphs of $G(\mathbb{Z}^d,C_\epsilon )$ and $G(\mathbb{Z}^d_{2n+1 + T},C_\epsilon )$ on vertex sets $[-n,n]^d$ are isomorphic.
	So, $\overline{\alpha}(G(\mathbb{Z}^d_{2n+1 + T} ,C_\epsilon )) \ge |I_n| / (2n+1+T)^d \ge  \frac{(2n+1)^d}{(2n+1+T)^d} (\overline{\alpha}(G(\mathbb{Z}^d ,C_\epsilon )) - \epsilon')$.
	Therefore, $\overline{\alpha}(G(\mathbb{Z}^d ,C_\epsilon ) \le \frac{(2n+1+T)^d}{(2n+1)^d}\overline{\alpha}(G(\mathbb{Z}^d_{2n+1 + T} ,C_\epsilon )) + \epsilon'$.
	Since $\epsilon' > 0$ can be chosen arbitrarily small and $\lim_{n \to \infty} \frac{(2n+1)^d}{(2n+1+T)^d} = 1$, we therefore have that  $\overline{\alpha}(G(\mathbb{Z}^d ,C_\epsilon ))   \le \limsup_{n \to \infty} \overline{\alpha}(G(\mathbb{Z}_n^d ,C_\epsilon ))$.
	As $G(\mathbb{Z}_n^d ,C_\epsilon)$ is a subgraph of $G(\mathbb{Z}_n^d ,C)$, we also have that $\overline{\alpha}(G(\mathbb{Z}^d ,C )) \le \overline{\alpha}(G(\mathbb{Z}^d ,C_\epsilon ))$.
	
	Finally, putting all of this together and examining limits, we get
	\begin{align*}
	\overline{\alpha}(G(\mathbb{Z}^d ,C ))
	& \le \lim_{\epsilon \to 0} \overline{\alpha}(G(\mathbb{Z}^d ,C_\epsilon ))   \\
	& \le \lim_{\epsilon \to 0} \limsup_{n \to \infty} \overline{\alpha}(G(\mathbb{Z}_n^d ,C_\epsilon ))   \\
	& \le \lim_{\epsilon \to 0} \frac{-\inf_{u\in (\mathbb{R} / \mathbb{Z})^d}  \widehat{w}_\epsilon(u) }{  \sup_{u\in (\mathbb{R} / \mathbb{Z})^d}  \widehat{w}_\epsilon(u) -\inf_{u\in (\mathbb{R} / \mathbb{Z})^d}  \widehat{w}_\epsilon(u)}   \\
	& = \frac{-\inf_{u\in (\mathbb{R} / \mathbb{Z})^d}  \widehat{w}(u) }{  \sup_{u\in (\mathbb{R} / \mathbb{Z})^d}  \widehat{w}(u)   -\inf_{u\in (\mathbb{R} / \mathbb{Z})^d}  \widehat{w}(u)},
	\end{align*}
	where the last line follows since for all $0 < \epsilon < s$ and $u\in (\mathbb{R} / \mathbb{Z})^d$, we have that $\widehat{w}(u) -\epsilon <   \widehat{w}_\epsilon(u) < \widehat{w}(u) + \epsilon$.
\end{proof}

\section{Distance graphs}

Theorem~\ref{p mod q} is trivial in the case that $p \equiv 0 \pmod{q}$, since $q\mathbb{Z} \times \{0\} \subset \mathbb{R}^2$ would be an infinite clique.
So, after scaling, it is enough to prove Theorem~\ref{p mod q} for coprime positive integers $p,q$ with $p<q$.
We start this section by choosing generating sets for our Cayley graphs of $\mathbb{Z}^2$.
The generating sets are chosen carefully to ensure that our Cayley graphs have large chromatic number, and that they still embed as a distance graph in the plane.

For each triple of positive integers $p,q,k$ with $p,q$ coprime and $p<q$, let
\[
D_{p,q,k}= \{\pm (q^{2k+t}-q^{2k-t} ,  -pq^{t} - 2q^{2k} ) : t=0,1,\ldots ,2k-1 \},
\]
and
\[
C_{p,q,k}= \bigcup_{j=0}^{\infty}  (jq+1) D_{p,q,k}  .
\]

The key property of the set $D_{p,q,k}$ that ensures that $G(\mathbb{Z}^2 , C_{p,q,k} )$ has large chromatic number is essentially the diversity of the highest power of $q$ that divides the two coordinates of the different points of $D_{p,q,k}$; as $t$ increases, the highest power of $q$ dividing the first coordinate decreases, and the highest power of $q$ dividing the second coordinate increases.
Intuitively, this is at least enough to avoid simple periodic colourings with a small period (depending on $k$).
We remark that if we defined $D_{p,q,k}$ more simply to be the set $\{\pm (q^{2k-t} ,  q^{t}) : t=0,1,\ldots ,2k-1 \}$, then $G(\mathbb{Z}^2 , C_{p,q,k} )$ would still have large chromatic number, however, then we would be unable to show that $G(\mathbb{Z}^2 , C_{p,q,k} )$ has a suitable embedding as a distance graph in the plane.

First, let us show that $G(\mathbb{Z}^2 , C_{p,q,k} )$ has an embedding as a distance graph in the plane.

\begin{lemma}\label{embed}
	For all positive integers $p,q,k$ with $p,q$ coprime and $p<q$, there is a injective function $h:\mathbb{Z}^2 \to \mathbb{R}^2$ such that $\|h(u)-h(v)\| \equiv p \pmod{q}$ whenever $u-v\in C_{p,q,k}$.
\end{lemma}

\begin{proof}
	Let $e_1=(1,0)$ and $e_2=\left(\frac{p}{2q^{2k}} , \sqrt{1 - \frac{p^2}{4q^{4k}}} \right)$.
	Let $h:\mathbb{Z}^2 \to \mathbb{R}^2$ be such that $h(x,y)=xe_1 + ye_2$.
	Clearly $h$ is injective, so it is enough to show that if $(x,y)$ is adjacent to $(x',y')$ in $G(\mathbb{Z}^2 , C_{p,q,k} )$, then the distance between $h(x,y)$ and $h(x',y')$ is equal to $p \pmod{q}$.
	
	So suppose that $(x,y)$ is adjacent to $(x',y')$ in $G(\mathbb{Z}^2,C_{p,q,k})$.
	Then $(x-x',y-y')\in C_{p,q,k}$.
	So there exists a non-negative integer $ j$ such that $(x-x',y-y')\in (jq + 1)D_{p,q,k}$.
	Therefore, there exists an integer $0\le t \le 2k-1$ such that $(x-x',y-y') = \pm (jq +1) (q^{2k+t}-q^{2k-t} ,  -pq^{t} - 2q^{2k} )$.
	Then
	\begin{align*}
	\frac{\| h(x,y) - h(x',y')  \|^2}{(jq+1)^2}
	& = \left[(q^{2k+t}-q^{2k-t}) + \frac{p}{2q^{2k}}  (-pq^{t} - 2q^{2k})  \right]^2 + \left[ \sqrt{1-\frac{p^2}{4q^{4k}}} (-pq^{t} - 2q^{2k}) \right]^2  \\
	& = (q^{2k+t}-q^{2k-t})^2 +  (-pq^{t} - 2q^{2k})^2 + \frac{p}{q^{2k}} (q^{2k+t}-q^{2k-t}) (-pq^{t} - 2q^{2k})  \\
	& = (q^{4k+2t}+ q^{4k-2t} -2q^{4k})  
	+ (p^2q^{2t} +4q^{4k} +4pq^{2k+t})  
	+ (-p^2q^{2t}  + 2pq^{2k-t}  -2 pq^{2k+t}  +p^2  )    \\
	& = q^{4k+2t} + q^{4k-2t}  + 2q^{4k}   +2pq^{2k+t} + 2pq^{2k-t}    +p^2 \\
	& = (q^{2k+t} +  q^{2k-t}   +p  )^2.
	\end{align*}
	Therefore $\|h(x,y)-h(x',y')\| = (jq+1)(q^{2k+t} +  q^{2k-t}   +p  )  \equiv p \pmod{q}$ as desired.
\end{proof}

Next, we will focus on applying Theorem~\ref{Z^d theta} to show that $G(\mathbb{Z}^2 , C_{p,q,k} )$ has large chromatic number for large $k$.
For our arguments, it is helpful to consider finite subsets of $C_{p,q,k}$.
For each quadruple of positive integers $p,q,k,n$ with $p,q$ coprime and $p<q$, let
\[
C_{p,q,k,n}= \bigcup_{j=0}^{n-1}  (jq+1) D_{p,q,k},
\]
and let $w_{p,q,k,n} : C_{p,q,k,n} \to \mathbb{R}$ be such that for each $c\in C_{p,q,k,n}$ of the form $(jq+1)d$ with $d\in D_{p,q,k}$, we have $w_{p,q,k,n}(c)=\frac{n-j}{n(n+1)}$.
As in Theorem \ref{Z^d theta}, for $u\in (\mathbb{R} / \mathbb{Z} )^2$, let
\[
\widehat{w}_{p,q,k,n}(u) = \sum_{x\in C_{p,q,k,n}} w_{p,q,k,n}(x) e^{2\pi i u \cdot x}.
\]
Clearly
\[
\sup_{u\in (\mathbb{R} / \mathbb{Z})^2}  \widehat{w}_{p,q,k,n}(u) = \widehat{w}_{p,q,k,n}(0) = \sum_{x\in C_{p,q,k,n}} w_{p,q,k,n}(x) = 4k\sum_{j=0}^{n-1} \frac{n-j}{n(n+1)} = 2k.
\]
So in order to effectively apply Theorem~\ref{Z^d theta}, it remains to bound $\inf_{u\in (\mathbb{R} / \mathbb{Z})^d}  \widehat{w}_{p,q,k,n}(u)$.

\begin{lemma}\label{bound}
	For all positive integers $p,q,k$ with $p,q$ coprime and $p<q$, we have
	\[
	\limsup_{{n} \to \infty} \inf_{u\in (\mathbb{R} / \mathbb{Z})^2}  \widehat{w}_{p,q,k,n}(u) \ge -2.
	\]
\end{lemma}

\begin{proof}
	Since $(\mathbb{R} / \mathbb{Z})^2$ is compact, for each positive integer $n$, there exists some $(a_n,b_n)\in (\mathbb{R} / \mathbb{Z})^2$ such that $\widehat{w}_{p,q,k,n}((a_n,b_n)) = \inf_{u\in (\mathbb{R} / \mathbb{Z})^2}  \widehat{w}_{p,q,k,n}(u)$.
	Furthermore, by the Bolzano-Weierstrass theorem, there is a strictly increasing sequence of positive integers $(n_{f})_{{f}=1}^{\infty}$ such that the sequence $\left( (a_{n_{f}} , b_{n_{f}})   \right)_{{f}=1}^\infty$ converges to some point $(a,b) \in (\mathbb{R} / \mathbb{Z})^2$.
	So for every $(x,y)\in D_{p,q,k}$, we have $\lim_{{f} \to \infty} a_{n_{f}}x + b_{n_{f}}y =ax+by$.
	
	Let,
	\begin{align*}
	& D_{p,q,k}^{\ref{lim f/q}} = \{(x,y) \in D_{p,q,k} : q(ax + by) \equiv 0 \pmod{1}, \ ax + by \not\equiv 0 \pmod{1} \},\\
	& D_{p,q,k}^{\ref{lim otherwise}} = \{(x,y) \in D_{p,q,k} : q(ax + by) \not\equiv 0 \pmod{1}\},\\
	& D_{p,q,k}^{\ref{lim 0}} = \{(x,y) \in D_{p,q,k} : ax + by \equiv 0 \pmod{1} \}.
	\end{align*}
	Observe that
	\begin{align*}
	\widehat{w}_{p,q,k,n_f}((a_{n_{f}},b_{n_{f}})) & = \sum_{(x,y)\in C_{p,q,k,n_{f}}} w((x,y)) e^{2\pi i (a_{n_{f}}x + b_{n_{f}}y) }\\
	& = \sum_{(x,y)\in D_{p,q,k}} \sum_{j=0}^{n_{f} -1} \frac{n_{f} -j}{n_{f} (n_{f} +1)} e^{2\pi i (jq+1) (a_{n_{f}}x + b_{n_{f}}y) }\\
	& = \sum_{(x,y)\in D_{p,q,k}} \frac{1}{n_{f} +1} \sum_{j=0}^{n_{f} -1} \frac{n_{f} -j}{n_{f} } e^{2\pi i (jq+1) (a_{n_{f}}x + b_{n_{f}}y) }\\
	& = \sum_{L\in \{1,2,3\}} \sum_{(x,y)\in D_{p,q,k}^{7.L}} \frac{1}{n_{f} +1} \sum_{j=0}^{n_{f} -1} \frac{n_{f} -j}{n_{f} } e^{2\pi i (jq+1) (a_{n_{f}}x + b_{n_{f}}y) }.
	\end{align*}
	We will obtain bounds for each of the summations over $D_{p,q,k}^{\ref{lim f/q}}, D_{p,q,k}^{\ref{lim otherwise}}, D_{p,q,k}^{\ref{lim 0}}$ individually. Note that each of $D_{p,q,k}^{\ref{lim f/q}}, D_{p,q,k}^{\ref{lim otherwise}}, D_{p,q,k}^{\ref{lim 0}}$ are centrally symmetric, so each of these three summations is real.

	We begin with bounding the summation over $D_{p,q,k}^{\ref{lim f/q}}$. Our careful choice of $D_{p,q,k}$ is what enables us to obtain such a bound for this summation.
	
	\begin{claim}\label{lim f/q}
		$\displaystyle
		\sum_{(x,y)\in D_{p,q,k}^{\ref{lim f/q}}} \frac{1}{n_{f} +1} \sum_{j=0}^{n_{f} -1} \frac{n_{f} -j}{n_{f} } e^{2\pi i (jq+1) (a_{n_{f}}x + b_{n_{f}}y) } \ge -2
		$.
	\end{claim}
	\begin{proof}
		Since $\frac{1}{n_{f} +1} \sum_{j=0}^{n_{f} -1} \frac{n_{f} -j}{n_{f} } =
		\frac{1}{2}$, it is enough to prove that $|D_{p,q,k}^{\ref{lim f/q}}|\le 4$.
		Suppose otherwise, then since $D_{p,q,k}^{\ref{lim f/q}}= - D_{p,q,k}^{\ref{lim f/q}}$, there exist integers $0\le  t_1 < t_2 < t_3     \le 2k-1$ such that for each $j\in \{1,2,3\}$, we have $T_j=(q^{2k+t_j}-q^{2k-t_j}, -pq^{t_j} - 2q^{2k}) \in D_{p,q,k}^{\ref{lim f/q}}$.
		Since neither of the three points $T_1,T_2,T_3$ is a scalar multiple of another, and $\mathbb{Z}^2$ does not contain $\mathbb{Z}^3$ as a subgroup, there exists non-zero integers $s_1,s_2,s_3$ with $\gcd(s_1,s_2,s_3)=1$, such that $s_1T_1 + s_2T_2 + s_3T_3 = (0,0)$.
		
		Observe that
		\begin{align*}
		0 & = s_1 (q^{2k+t_1}-q^{2k-t_1}) +
		s_2 (q^{2k+t_2}-q^{2k-t_2}) +
		s_3 (q^{2k+t_3}-q^{2k-t_3})\\
		& = s_1 q^{2k-t_1} (q^{2t_1}-1) +
		s_2 q^{2k-t_2} (q^{2t_2}-1) +
		s_3 q^{2k-t_3} (q^{2t_3}-1).
		\end{align*}
		Since $q$ and $q^{2t_j}-1$ are coprime for each $j\in\{1,2,3\}$, and $2k-t_3 < 2k-t_2 < 2k-t_1$, it follows that $q$ divides $s_3$.
		Similarly,
		\begin{align*}
		0 & = s_1 (-pq^{t_1} - 2q^{2k}) +
		s_2 (-pq^{t_2} - 2q^{2k}) +
		s_3 (-pq^{t_3} - 2q^{2k})\\
		& = s_1 q^{t_1} (-p - 2q^{2k-t_1}) +
		s_2 q^{t_2} (-p - 2q^{2k-t_2}) +
		s_3 q^{t_3} (-p - 2q^{2k-t_3}).
		\end{align*}
		So $q$ divides $s_1$ since $q$ and $p + 2q^{2k-t_j}$ are coprime for each $j\in\{1,2,3\}$, and $t_1<t_2<t_3$.
		As $\gcd(s_1,s_2,s_3)=1$, it now follows that $q$ and $s_2$ are coprime.

		By definition of $D_{p,q,k}^{\ref{lim f/q}}$, for each $j\in \{1,2,3\}$, we have $qT_j \cdot (a,b) \equiv 0 \pmod{1}$.
		So, as $q$ divides both $s_1$ and $s_3$, we have that $s_1T_1 \cdot (a,b) \equiv 0 \pmod{1}$, and $s_3T_3 \cdot (a,b) \equiv 0 \pmod{1}$.
		Since $s_1T_1 + s_2T_2 + s_3T_3 = (0,0)$, this implies that $s_2T_2 \cdot (a,b) \equiv 0 \pmod{1}$. So $q ( T_2 \cdot (a,b) ) \equiv s_2 ( T_2 \cdot (a,b) ) \equiv 0 \pmod{1}$.
		As $q$ and $s_2$ are coprime, this implies that $T_2 \cdot (a,b) \equiv 0 \pmod{1}$, contracting the definition of $D_{p,q,k}^{\ref{lim f/q}}$.
		Thus $|D_{p,q,k}^{\ref{lim f/q}}| \le 4$ as required.
	\end{proof}
	
	Next we tackle the summation over $D_{p,q,k}^{\ref{lim otherwise}}$.
	It is simply the natural choice of $C_{p,q,k,n}$ given $D_{p,q,k}$ that enables us to obtain such an estimate for this summation.
	
	\begin{claim}\label{lim otherwise}
		$\displaystyle
		\lim_{f \to \infty} \sum_{(x,y)\in D_{p,q,k}^{\ref{lim otherwise}}} \frac{1}{n_{f} +1} \sum_{j=0}^{n_{f} -1} \frac{n_{f} -j}{n_{f} } e^{2\pi i (jq+1) (a_{n_{f}}x + b_{n_{f}}y) } 
		= 0$.
	\end{claim}
	\begin{proof}
		Since $D_{p,q,k}^{\ref{lim otherwise}} \subseteq D_{p,q,k}$ is a finite set,
		it is enough to show that if $(x,y)\in D_{p,q,k}^{\ref{lim otherwise}}$, then
		\[
		\lim_{f \to \infty} \frac{1}{n_{f} +1} \sum_{j=0}^{n_{f} -1} \frac{n_{f} -j}{n_{f} } e^{2\pi i (jq+1) (a_{n_{f}}x + b_{n_{f}}y) } 
		= 0.
		\]
		
		Choose a positive integer $m$ arbitrarily.
		For each $s\in \{1,\ldots, m\}$, let $M_s= \mathbb{N} \cap \big[ \frac{(s-1)n_{f} }{m} , \frac{sn_{f} }{m} \big)$.
		Observe that
		\begin{align*}
		& \left| \frac{1}{n_{f} +1} \sum_{j=0}^{n_{f} -1} \frac{n_{f} -j}{n_{f} } e^{2\pi i (jq+1) (a_{n_{f}}x + b_{n_{f}}y) } \right|\\
		& \ \ \ \ \ \ \ \ \le 	\frac{1}{n_{f}}  \sum_{s=1}^m \left| \sum_{j\in M_s} \frac{n_{f} -j}{n_{f} } e^{2\pi i jq (a_{n_{f}}x + b_{n_{f}}y) } \right| \\
		& \ \ \ \ \ \ \ \ = \frac{1}{n_{f} }  \sum_{s=1}^m \left|  \sum_{j\in M_s} \frac{\max\{M_s\} - j}{n_{f} } e^{2\pi i jq+ (a_{n_{f}}x + b_{n_{f}}y) }
		+ \sum_{j\in M_s} \frac{n_{f} -\max\{M_s\}}{n_{f} } e^{2\pi i jq (a_{n_{f}}x + b_{n_{f}}y) }   \right|\\
		& \ \ \ \ \ \ \ \ \le \frac{1}{n_{f} }  \sum_{s=1}^m \left( \left|  \sum_{j\in M_s} \frac{\max\{M_s\} - j}{n_{f} } e^{2\pi i jq (a_{n_{f}}x + b_{n_{f}}y) } \right|
		+   \left| \sum_{j\in M_s} \frac{n_{f} -\max\{M_s\}}{n_{f} } e^{2\pi i jq (a_{n_{f}}x + b_{n_{f}}y) }   \right| \right)\\
		& \ \ \ \ \ \ \ \ \le \frac{1}{n_{f} }  \sum_{s=1}^m    \sum_{j\in M_s} \frac{\max\{M_s\} - \min\{M_s\}}{n_{f} } 
		+  \frac{1}{n_{f} }  \sum_{s=1}^m  \left| \sum_{j\in M_s} e^{2\pi i jq (a_{n_{f}}x + b_{n_{f}}y) }   \right|\\
		& \ \ \ \ \ \ \ \ = \frac{1}{n_{f} }  \sum_{s=1}^m    \sum_{j\in M_s} \frac{|M_s|}{n_{f} } 
		+  \frac{1}{n_{f} }  \sum_{s=1}^m  \left| \sum_{j\in M_s} e^{2\pi i (j -\min\{M_s\} )q (a_{n_{f}}x + b_{n_{f}}y) }   \right|\\
		& \ \ \ \ \ \ \ \ = \frac{1}{n_{f}^2 }  \sum_{s=1}^m     |M_s|^2 
		+  \frac{1}{n_{f} }  \sum_{s=1}^m  \frac{\left| 1 - e^{2\pi i |M_s| q (a_{n_{f}}x + b_{n_{f}}y) } \right| }{ \left| 1 - e^{2\pi i q (a_{n_{f}}x + b_{n_{f}}y) } \right| } \\
		& \ \ \ \ \ \ \ \ \le \frac{m}{n_{f}^2 }  \left(\frac{n_{f}}{m}+1 \right)^2
		+  \frac{m}{n_{f} }  \left( \frac{ 1 }{ \left| 1 - e^{2\pi i q (a_{n_{f}}x + b_{n_{f}}y) } \right|  } \right)
		\end{align*}
		
		By definition of $D_{p,q,k}^{\ref{lim otherwise}}$, we have $q(ax+by) \not\equiv 0 \pmod{1}$. So there exists some $\epsilon > 0$ and $F\in \mathbb{N}$, such that $q(a_{n_{f}}x + b_{n_{f}}y) \in [\epsilon , 1 - \epsilon] \pmod{1}$ for all ${f} \ge F$.
		Since $ \left|  1 - e^{2\pi i q (a_{n_{f}}x + b_{n_{f}}y) }       \right| $ is bounded away from $0$ for ${f} \ge F$, we get that
		\[
		\lim_{{f} \to \infty} \left| \frac{1}{n_{f} +1} \sum_{j=0}^{n_{f} -1} \frac{n_{f} -j}{n_{f} } e^{2\pi i (jq+1) (a_{n_{f}}x + b_{n_{f}}y) } \right| 
		 \le  \lim_{{f} \to \infty} \frac{m}{n_{f}^2 }  \left(\frac{n_{f}}{m}+1 \right)^2
		+  \lim_{{f} \to \infty} \frac{m}{n_{f} }  \left( \frac{ 1 }{ \left| 1 - e^{2\pi i q (a_{n_{f}}x + b_{n_{f}}y) } \right|  } \right)
		 =  \frac{1}{m}.
		\]

		Since $m$ was chosen arbitrarily, the claim now follows.
	\end{proof}
	
	Lastly, we handle the summation over $D_{p,q,k}^{\ref{lim 0}}$. Here, the choice of our function $w_{p,q,k,n}$, given $C_{p,q,k,n}$ is what enables the following bound.

	\begin{claim}\label{lim 0}
		$\displaystyle
		\liminf_{f \to \infty} \sum_{(x,y)\in D_{p,q,k}^{\ref{lim 0}}} \frac{1}{n_{f} +1} \sum_{j=0}^{n_{f} -1} \frac{n_{f} -j}{n_{f} } e^{2\pi i (jq+1) (a_{n_{f}}x + b_{n_{f}}y) } 
		\ge 0$.
	\end{claim}
	\begin{proof}
		Since $D_{p,q,k}^{\ref{lim 0}}\subseteq D_{p,q,k}$ is a finite centrally symmetric set, by the identity $e^{i\theta} + e^{-i\theta}=2\cos(\theta)$, it is enough to show that if $(x,y) \in D_{p,q,k}^{\ref{lim 0}}$, then 
		\[
		\liminf_{f \to \infty} \frac{1}{n_{{f}} +1 } \sum_{j=0}^{n_{f} -1} \frac{n_{f} -j}{n_{f} } \cos(2\pi (jq+1) (a_{n_{f}}x + b_{n_{f}}y) ) 
		\ge 0.
		\]
		For each ${f} \ge 1$, let $\theta_{f} \in [0,1/2]$ be such that $a_{n_{f}}x + b_{n_{f}}y  \equiv \pm \theta_{f} \pmod{1}$. So $\lim_{{f} \to \infty} \theta_{f} = 0$ by the definition of $D_{p,q,k}^{\ref{lim 0}}$.
		If $\theta_{f}=0$, then
		\[
		\frac{1}{n_{{f}} +1 } \sum_{j=0}^{n_{f} -1} \frac{n_{f} -j}{n_{f} } \cos(2\pi (jq+1) \theta_{f} ) 
		= \frac{1}{n_{{f}} +1 } \sum_{j=0}^{n_{f} -1} \frac{n_{f} -j}{n_{f} }
		= \frac{1}{2} > 0.
		\]
		So, we may assume that $\theta_{f} \not= 0$.

		Observe that
		\[
		\lim_{{f} \to \infty}  \left( 
			\frac{1}{n_{{f}} +1 } \sum_{j=0}^{n_{f} -1} \frac{n_{f} -j}{n_{f} } \cos(2\pi (jq+1) \theta_{f} ) 
		-
			\frac{1}{ n_{{f}}q \theta_{f} } \int_{0}^{ n_{{f}}q \theta_{f}  }
			\left( 1 - \frac{t}{  n_{{f}}q \theta_{f} } \right) \cos(2 \pi  t) \,dt
		\right)
		=0.
		\]
		Therefore,
		\[
		\liminf_{f \to \infty} \frac{1}{n_{{f}} +1 } \sum_{j=0}^{n_{f} -1} \frac{n_{f} -j}{n_{f} } \cos(2\pi (jq+1) \theta_{f} )
		= 
		\liminf_{f \to \infty} \frac{1}{ n_{{f}}q \theta_{f} } \int_{0}^{ n_{{f}}q \theta_{f}  }
		\left( 1 - \frac{t}{  n_{{f}}q \theta_{f} } \right) \cos(2 \pi  t) \,dt.
		\]
		Let us evaluate the above integral, we substitute $\ell= 2\pi n_{{f}}q \theta_{f}$ for convenience (note that $\ell\not=0$),
		\begin{align*}
		\frac{2\pi}{  \ell  }   \int_{0}^{ \frac{\ell}{ 2\pi }  }
		\left( 1 - \frac{2\pi t}{  \ell  } \right) \cos(2\pi t) \,dt 
		& =\frac{1}{  \ell  }   \int_{0}^{ \ell  }
		\left( 1 - \frac{t}{  \ell  } \right) \cos(t) \,dt \\
		& = \frac{1}{  \ell  }   \int_{0}^{ \ell   } \cos(t) \,dt
		- \frac{1}{ \ell^2 }   \int_{0}^{ \ell   }
		t \cos(t) \,dt\\
		& = \frac{ \sin(  \ell  )  }{  \ell  }
		- \frac{ \sin(  \ell  )  }{  \ell  }
		+ \frac{1}{  \ell^2  } \int_{0}^{ \ell   } \sin(t) \,dt\\
		& = \frac{1 - \cos( \ell )  }{  \ell^2  }.
		\end{align*}
		The claim now follows since $\frac{1 - \cos( \ell )  }{  \ell^2  } \ge 0$.
	\end{proof}

	By Claim \ref{lim f/q}, Claim \ref{lim otherwise}, and Claim \ref{lim 0}, 
	\[
	\limsup_{{n} \to \infty} \inf_{u\in (\mathbb{R} / \mathbb{Z})^2}  \widehat{w}_{p,q,k,n}(u) \ge \liminf_{{f} \to \infty} \inf_{u\in (\mathbb{R} / \mathbb{Z})^2}  \widehat{w}_{p,q,k,n_f}(u) \ge -2,
	\]
	as desired.
\end{proof}

We are now ready to prove that $G(\mathbb{Z}^2, C_{p,q,k} )$ has large chromatic number (when $k$ is large).

\begin{theorem}\label{graphs}
	For all positive integers $p,q,k$ with $p,q$ coprime and $p<q$, we have
	$\overline{\alpha}(G(\mathbb{Z}^2, C_{p,q,k} )  \le \frac{1}{k+1}$.
	As a consequence, $\chi( G(\mathbb{Z}^2, C_{p,q,k} )  ) \ge k+1$.
\end{theorem}

\begin{proof}
	As discussed before Lemma~\ref{bound}, we have $\sup_{u\in (\mathbb{R} / \mathbb{Z})^2}  \widehat{w}_{p,q,k,n}(u) = 2k$ for every positive integer $n$.
	By Lemma~\ref{bound}, we have $\limsup_{{n} \to \infty} \inf_{u\in (\mathbb{R} / \mathbb{Z})^2}  \widehat{w}_{p,q,k,n}(u) \ge -2$.
	Observe that for every positive integer $n$, we have that $\inf_{u\in (\mathbb{R} / \mathbb{Z})^2}  \widehat{w}_{p,q,k,n}(u)\le 0$, since 
	$\int_{[0,1]^2} \widehat{w}_{p,q,k,n}(u) \,du = 0$.
	Therefore, by Theorem~\ref{Z^d theta}, we have
	\begin{align*}
	\overline{\alpha}(G(\mathbb{Z}^2, C_{p,q,k} )  
	& \le \liminf_{{n} \to \infty} \overline{\alpha}(G(\mathbb{Z}^2, C_{p,q,k,n} )   \\
	& \le \liminf_{{n} \to \infty} 
	\frac{-\inf_{u\in (\mathbb{R} / \mathbb{Z})^2}  \widehat{w}_{p,q,k,n}(u) }{  \sup_{u\in (\mathbb{R} / \mathbb{Z})^2}  \widehat{w}_{p,q,k,n}(u)   -\inf_{u\in (\mathbb{R} / \mathbb{Z})^2}  \widehat{w}_{p,q,k,n}(u)},
	\\
	& =  
	\frac{-\limsup_{{n} \to \infty} \inf_{u\in (\mathbb{R} / \mathbb{Z})^2}  \widehat{w}_{p,q,k,n}(u) }{  2k   -\limsup_{{n} \to \infty} \inf_{u\in (\mathbb{R} / \mathbb{Z})^2}  \widehat{w}_{p,q,k,n}(u)},
	\\
	& \le \frac{1}{k+1}.
	\end{align*}
	As a consequence, $\chi( G(\mathbb{Z}^2, C_{p,q,k} )  ) \ge 1/\overline{\alpha}(G(\mathbb{Z}^2, C_{p,q,k} ) \ge k+1$.
\end{proof}

As previously discussed at the start of the section, Theorem~\ref{p mod q} is trivial in the case that $p \equiv 0 \pmod{q}$, since $q\mathbb{Z} \times \{0\} \subset \mathbb{R}^2$ would be an infinite clique.
By scaling, we may then reduce to the case that $p$ and $q$ are coprime.
So, Theorem~\ref{p mod q} now follows from Theorem~\ref{graphs} and Lemma~\ref{embed}.

We finish the paper with a proof that the graphs in Theorem~\ref{graphs} are triangle-free, and thus provide a new construction of triangle-free graphs with arbitrarily large chromatic number.

\begin{proposition}\label{triangle-free}
	For all positive integers $p,q,k$ with $p,q$ coprime and $p<q$, the graph $G(\mathbb{Z}^2, C_{p,q,k})$ is triangle-free.
\end{proposition}

\begin{proof}
	Suppose otherwise, then there exists $t_1,t_2,t_3\in \{0,1,\ldots , 2k-1\}$, and non-negative integers $j_1,j_2,j_3$ such that $(j_1q+1)T_1=(j_2q+1)T_2+(j_3q+1)T_3$, where $T_i=(q^{2k+t_i}-q^{2k-t_i}, -pq^{t_i} - 2q^{2k})$ for each $i\in \{1,2,3\}$. Since $p\not\equiv 2p \pmod{q}$, it follows that $t_1,t_2,t_3$ are not all equal. Therefore, $T_1,T_2,T_3,(0,0)$ are not collinear. Since $(j_1q+1)T_1=(j_2q+1)T_2+(j_3q+1)T_3$, it follows that no two of $T_1,T_2,T_3$ are collinear with $(0,0)$. So $t_1,t_2,t_3$ are distinct, and without loss of generality, we may assume that $t_3>t_2$. The highest power of $q$ dividing $(j_1q+1)(-pq^{t_1} - 2q^{2k})$ is $q^{t_1}$, while the highest power of $q$ dividing $(j_2q+1)(-pq^{t_2} - 2q^{2k}) + (j_3q+1)(-pq^{t_3} - 2q^{2k})$ is $q^{t_2}$. This contradicts the fact that $(j_1q+1)T_1=(j_2q+1)T_2+(j_3q+1)T_3$.
\end{proof}

\section*{Acknowledgements}

The author thanks Sabrina Lato, Rose McCarty, and Micha\l{} Pilipczuk for helpful discussions on Rosenfeld's odd distance problem and its solution. The author also thanks Chris Godsil for bringing~\cite{lovasz1975spectra} to their attention, and the anonymous referees for valuable comments that have improved the paper.

\bibliographystyle{amsplain}

\begin{thebibliography}{10}
	
	\bibitem{ardal2009odd}
	Hayri Ardal, J{\'a}n Ma{\v{n}}uch, Moshe Rosenfeld, Saharon Shelah, and
	Ladislav Stacho, \emph{The odd-distance plane graph}, Discrete \&
	Computational Geometry \textbf{42} (2009), no.~2, 132--141.
	
	\bibitem{bachoc2014spectral}
	Christine Bachoc, Evan DeCorte, Fernando~M{\'a}rio de~Oliveira~Filho, and Frank
	Vallentin, \emph{Spectral bounds for the independence ratio and the chromatic
		number of an operator}, Israel Journal of Mathematics \textbf{202} (2014),
	no.~1, 227--254.
	
	\bibitem{beineke2015topics}
	Lowell~W Beineke and Robin~J Wilson, \emph{Topics in {C}hromatic {G}raph
		{T}heory}, vol. 156, Cambridge University Press, 2015.
	
	\bibitem{bourgain1986szemeredi}
	Jean Bourgain, \emph{A {S}zemer{\'e}di type theorem for sets of positive
		density in $\mathbb{R}^k$}, Israel Journal of Mathematics \textbf{54} (1986),
	no.~3, 307--316.
	
	\bibitem{chybowska2022coloring}
	Joanna Chybowska-Sok{\'o}{\l}, Konstanty Junosza-Szaniawski, and Krzysztof
	W{{e}}sek, \emph{Coloring distance graphs on the plane}, Discrete Mathematics \textbf{346} (2023), no. 7, 113441.
	
	\bibitem{damasdi2021odd}
	G{\'a}bor Dam{\'a}sdi, \emph{Odd wheels are not odd-distance graphs}, Discrete
	\& Computational Geometry (2021), 1--11.
	
	\bibitem{davies2023prime}
	James Davies, Rose McCarty, and Micha{\l} Pilipczuk, \emph{Prime and polynomial distances in colourings of the plane}, arXiv preprint
	arXiv:2308.02483 (2023).
	
	\bibitem{de2018chromatic}
	Aubrey D. N.~J. de~Grey, \emph{The chromatic number of the plane is at least
		5}, Geombinatorics \textbf{28} (2018), no.~1, 18--31.
	
	\bibitem{de2010fourier}
	Fernando~M{\'a}rio de~Oliveira~Filho and Frank Vallentin, \emph{Fourier
		analysis, linear programming, and densities of distance avoiding sets in
		$\mathbb{R}^n$}, Journal of the European Mathematical Society \textbf{12}
	(2010), no.~6, 1417--1428.
	
	\bibitem{dutour2021coloring}
	Mathieu Dutour~Sikiri{\'c}, David~A Madore, Philippe Moustrou, and Frank
	Vallentin, \emph{Coloring the {V}oronoi tessellation of lattices}, Journal of
	the London Mathematical Society \textbf{104} (2021), no.~3, 1135--1171.
	
	\bibitem{erdos1994twenty}
	Paul Erd\H{o}s, \emph{Twenty five years of questions and answers}, 25th
	Southeastern International Conference on Combinatorics, Graph Theory and
	Computing, Boca Raton, Florida, 1994.
	
	\bibitem{falconer1986plane}
	Kenneth~J Falconer and John~M Marstrand, \emph{Plane sets with positive density
		at infinity contain all large distances}, Bulletin of the London Mathematical
	Society \textbf{18} (1986), no.~5, 471--474.
	
	\bibitem{furstenberg1990ergodic}
	Hillel F{\"u}rstenberg, Yitzchak Katznelson, and Benjamin Weiss, \emph{Ergodic
		theory and configurations in sets of positive density}, Mathematics of Ramsey
	theory, Springer, 1990, pp.~184--198.
	
	\bibitem{golovanov2023odd}
	Alexander Golovanov, Andrey Kupavskii, and Arsenii Sagdeev, \emph{Odd-distance
		and right-equidistant sets in the maximum and {M}anhattan metrics}, European
	Journal of Combinatorics \textbf{107} (2023), 103603.
	
	\bibitem{graham2010open}
	Ron Graham and Eric Tressler, \emph{Open problems in {E}uclidean {R}amsey
		{T}heory}, Ramsey {T}heory: {Y}esterday, {T}oday, and {T}omorrow \textbf{285}
	(2010), 115--120.
	
	\bibitem{heuleeasier}
	Marijn J.~H. Heule, \emph{Easier variants of notorious math problems}, Nieuw
	Archief voor Wiskunde \textbf{5/22} (2021), no.~3, 153--157.
	
	\bibitem{heuleodd}
	Marijn J.~H. Heule, \emph{Odd-distance virtual edges in unit-distance graphs},
	Geombinatorics \textbf{31} (2021), no.~2, 68--76.
	
	\bibitem{jensen2011graph}
	Tommy~R Jensen and Bjarne Toft, \emph{Graph {C}oloring {P}roblems}, John Wiley
	\& Sons, 2011.
	
	\bibitem{kalai15some}
	Gil Kalai, \emph{Some old and new problems in combinatorial geometry {I}:
		{A}round {B}orsuk's problem}, Surveys in Combinatorics (Artur Czumaj, Agelos
	Georgakopoulos, Daniel Kr\'a{\v{l}}, Vadim Lozin, and Oleg Pikhurko, eds.),
	London Mathematical Society Lecture Note Series, vol. 424, Cambridge
	University Press, 2015, pp.~147--174.
	
	\bibitem{kechris2016descriptive}
	Alexander~S Kechris and Andrew~S Marks, \emph{Descriptive graph combinatorics},
	(Preprint) (2016).
	
	\bibitem{lovasz1975spectra}
	L{\'a}szl{\'o} Lov{\'a}sz, \emph{Spectra of graphs with transitive groups},
	Periodica Mathematica Hungarica \textbf{6} (1975), no.~2, 191--195.
	
	\bibitem{lovasz1979shannon}
	L{\'a}szl{\'o} Lov{\'a}sz, \emph{On the {S}hannon capacity of a graph}, IEEE Transactions on
	Information theory \textbf{25} (1979), no.~1, 1--7.
	
	\bibitem{parts20226}
	Jaan Parts, \emph{A 6-chromatic odd-distance graph in the plane},
	Geombinatorics \textbf{31} (2022), no.~3, 124--137.
	
	\bibitem{quas2009distances}
	Anthony Quas, \emph{Distances in positive density sets in $\mathbb{R}^n$},
	Journal of Combinatorial Theory, Series A \textbf{116} (2009), no.~4,
	979--987.
	
	\bibitem{rosenfeld1996odd}
	Moshe Rosenfeld, \emph{Odd integral distances among points in the plane},
	Geombinatorics \textbf{5} (1996), no.~4, 156--159.
	
	\bibitem{rosenfeld2008some}
	Moshe Rosenfeld, \emph{Some of my favorite ``lesser known'' problems}, Ars Mathematica
	Contemporanea \textbf{1} (2008), no.~2, 137--143.
	
	\bibitem{rosenfeld2014forbidden}
	Moshe Rosenfeld and Nam~L{\^e} Ti{\^e}n, \emph{Forbidden subgraphs of the
		odd-distance graph}, Journal of Graph Theory \textbf{75} (2014), no.~4,
	323--330.
	
	\bibitem{scott2020survey}
	Alex Scott and Paul Seymour, \emph{A survey of $\chi$-boundedness}, Journal of
	Graph Theory \textbf{95} (2020), no.~3, 473--504.
	
	\bibitem{soifer2010between}
	Alexander Soifer, \emph{Between the line and the plane: Chromatic {\'e}tude in
		6 movements.}, Mathematics Competitions \textbf{23} (2010), no.~2, 30--45.
	
	\bibitem{soifer2016hadwiger}
	Alexander Soifer, \emph{The {H}adwiger--{N}elson problem}, Open {P}roblems in
	{M}athematics, Springer, 2016, pp.~439--457.
	
	\bibitem{steinhardt2009coloring}
	Jacob Steinhardt, \emph{On coloring the odd-distance graph}, The Electronic
	Journal of Combinatorics \textbf{16} (2009), no.~1, N12.
	
\end{thebibliography}

\end{document}